\documentclass[12pt,a4paper,reqno]{amsart}
\usepackage{amssymb,amsmath}
\usepackage{hyperref}
\usepackage{setspace}
\usepackage{graphicx}
\usepackage{amsmath, array, makecell}
\usepackage{cellspace}

\theoremstyle{definition}

\begin{document}

\begin{center}
{\Large Conservation Laws And The Applicability Of Group Theoretical Technique to Non-Linear  Chaffee-Infante Equation}

\bigskip

\textbf{Muhammad Naqeeb}

Department of Mathematics, Quaid-i-Azam University, Islamabad-44000, Pakistan.

Email; mnaqeem@math.qau.edu.pk

\bigskip
\bigskip

\textbf{ Amjad Hussain}

Department of Mathematics, Quaid-i-Azam University, Islamabad-44000, Pakistan.

Email; a.hussain@qau.edu.pk

\bigskip
\bigskip

\end{center}

\begin{center}
{\small Abstract}
\end{center}
   Group theoretical technique is applied for the sake of similarity analysis and conservation laws of non-linear Chaffee-Infante equation. The similarity reductions are performed to get the exact analytic solutions by power series method, through the  construction of its infinitesimal generators and optimal system.   \par
keywords: Group-theoretical method, Chaffee-Infante equation, conservation laws, Power series, Invariant and transformed solutions.
\section{introduction}
Initiated by Sophus Lie, The group theoretical analysis has deep roots for handling linear and highly non-linear differential equations\cite{NNN}, specifically for exact solutions. From Lie-point symmetry to its extended forms like approximate symmetry, Backlund or contact symmetry , generalized symmetry and hidden symmetry. It is widely applicable on recent advance types of stochastic equations, reaction diffusion equations and  fractional equations. It explores not only solutions of different kinds, but opens variety of fields to step in from conservation laws of dynamics \cite{EMN}, differential manifolds of differential geometry, Models of Astrophysics\cite{ACTTTTTT}, to symmetries of differential equations\cite{DJAF}.\par
In this paper, analyzing non-linear Chaffee-Infante equation,
\begin{equation}
u_t-u_{xx}+\lambda(u^3-u)=0, \quad \lambda\ne 0.
\end{equation}
by group theoretical method, which is one of the reaction diffusion equations modeled in Mathematical physics\cite{YM}, and occurs as mathematical model to several physical phenomenon. we have understood conservation laws and various solutions, which are invariant, transformed, fundamental and solutions generated by infinite ideal\cite{IB}. by grasping the concept of optimal system\cite{BL} and properties of lie algebra\cite{PH} of non-linear Chaffee-Infante equation.  The detailed discussion is done on corresponding solutions of infinitesimal groups generated by the fundamental theorem of linear combinations of infinitesimal generators. The possibility of constructing conservation laws associated with symmetries was first performed by Emmy Noether for Euler-Lagrange equations. The concept was extended by Ibragimov for partial differential equations. For Chaffee-Infante equation Conservation laws\cite{Ibragimov} are constructed by analyzing its adjointness, self-adjointness, quasi-self adjointness\cite{DRAMA}, weak self-adjointness\cite{DRAMAM} and non-linear adjointness.
\section{Preliminaries}
Lie symmetry analysis being a powerful group theoretical technique\cite{PJO} has been employed to obtain the infinitesimal generator and commutation table of lie algebra for the Chaffee-Infante equation. The one-parameter infinitesimal transformation\cite{THES} are given as;
\begin{equation*}\begin{aligned}&
 x \rightarrow x+\epsilon\omega(x,t,u), \\&
   t \rightarrow t+\epsilon\psi(x,t,u),\\&
   u \rightarrow u+\epsilon\chi(x,t,u).  \end{aligned}
 \end{equation*}
 The infinitesimal operator or group generator associated with above transformations is;
 \begin{equation}
 G=\omega(x,t,u)\frac{\partial }{\partial x}+ \psi(x,t,u)\frac{\partial }{\partial t}+\chi(x,t,u)\frac{\partial }{\partial u}.
 \end{equation}
 As of our equation, is second order. So the second prolongation of infinitesimal generator would be,
 \begin{equation}\begin{aligned}&
 G^{(2)}=\omega(x,t,u)\frac{\partial }{\partial x}+ \psi(x,t,u)\frac{\partial }{\partial t}+\chi(x,t,u)\frac{\partial }{\partial u}+\chi_{t}(x,t,u)\frac{\partial }{\partial u_t}\\&+\chi_{x}(x,t,u)\frac{\partial }{\partial u_x}
 +\chi_{xx}(x,t,u)\frac{\partial }{\partial u_{xx}}+\chi_{tx}(x,t,u)\frac{\partial }{\partial u_{tx}}+\chi_{tt}(x,t,u)\frac{\partial }{\partial u_{tt}}.\end{aligned}
 \end{equation}
whereas the invariance condition of the symmetry is given as,
\begin{equation}
G^{(2)}\rho=0,  \quad \rho=u_t-u_{xx}+\lambda (u^3-u).
\end{equation}
The extended transformations are given as;
\begin{equation}\begin{aligned}&
\chi_{t}=D_{t}(\chi)-u_{t}D_{t}(\psi)-u_{x}D_{t}(\omega), \\&
\chi_{x}=D_{x}(\chi)-u_{t}D_{x}(\psi)-u_{x}D_{x}(\omega),  \\&
\chi_{tt}=D_{t}(\chi_t)-u_{tt}D_{t}(\psi)-u_{tx}D_{t}(\omega),  \\&
\chi_{xx}=D_{x}(\chi_x)-u_{tx}D_{x}(\psi)-u_{xx}D_{x}(\omega),  \\&
\chi_{tx}=D_{x}(\chi_t)-u_{tt}D_{x}(\psi)-u_{tx}D_{x}(\omega). \end{aligned}
\end{equation}
Whereas, $D_x$ and $D_t$ are given as,
\begin{equation}
D_{t}=\frac{\partial}{\partial t}+u_t\frac{\partial}{\partial u}+u_{tt}\frac{\partial}{\partial u_t}+u_{tx}\frac{\partial}{\partial u_x}+...,
\end{equation}
\begin{equation}
D_{x}=\frac{\partial}{\partial x}+u_x\frac{\partial}{\partial u}+u_{xx}\frac{\partial}{\partial u_x}+u_{tx}\frac{\partial}{\partial u_t}+...
\end{equation}

\section{Lie symmetries of Chaffee-Infante equation}
To calculate Lie point symmetries corresponding to infinitesimal generator. Employ invariance condition would lead to determining equations, which are computed as;
\begin{equation}
\chi_t-\chi_{xx}+(\lambda u^{3}-u)(-\chi_u+2\omega_x)+(3u^{2}\lambda-\lambda)\chi=0,
\end{equation}
\begin{equation}
\psi_{xx}-\psi_{t}+2\omega_{x}+(\lambda u^{3}-\lambda u)\psi_u=0,
\end{equation}
\begin{equation}
\omega_{xx}-\omega_{t}-2\chi_{xu}+3 \omega_{u}(\lambda u^{3}-\lambda u)=0,
\end{equation}
\begin{equation}
\psi_{xu}+\omega_{u}=0,
\end{equation}
\begin{equation}
\omega_{uu}=0,
\end{equation}
\begin{equation}
\psi_{uu}=0,
\end{equation}
\begin{equation}
\psi_{u}=0,
\end{equation}
\begin{equation}
\psi_{x}=0,
\end{equation}
\begin{equation}
2\psi_{xu}-\chi_{uu}=0.
\end{equation}
Simplifying the above system gives our 'system of determining equations' in simplest form which is;
\begin{equation}\begin{aligned}&
\psi_{x}=0, \quad \psi_{u}=0, \quad \omega_{u}=0, \quad \chi_{uu}=0, \\&
\chi_{t}-\chi_{xx}-(\lambda u^{3}-u)(\chi_{u}-2\omega_{x})+(3u^{2}\lambda-\lambda)\chi=0, \\&
2\omega_{x}-\psi_{t}=0, \quad \omega_{xx}-\omega_{t}-2\chi_{xu}=0. \end{aligned}
\end{equation}
Solving the above system for the infinitesimals that would lead to find the symmetries and infinitesimal generators of the Chaffee-Infantee equation, These are the required infinitesimals of the symmetries to be considered so far;
\begin{equation}
\psi=A(t),
\end{equation}
\begin{equation}
\omega=B(t,x),
\end{equation}
\begin{equation}
\chi=P(t,x)u+Q(t,x).
\end{equation}
Which further simplifies the above system given as;
\begin{equation}\begin{aligned}&
2B_{x}-A_{t}=0,\quad B_{xx}-B_{t}-2P_{x}=0,  \\&
(\lambda u^{3}-u)(p-2B_{x})-(3u^{2}\lambda-\lambda)(Pu+Q)=(P_{t}-P_{xx})u+(Q_{t}-Q_{xx}). \end{aligned}
\end{equation}
where A, B, Q, P are arbitrary functions.
To compute infinitesimals we could consider various cases, corresponding to above mentioned constant functions of independent variables,
considering P=0,  Q=0, the corresponding infinitesimals are;
\begin{equation}
\psi=C_{1},\quad \omega=C_{2},\quad \chi=0.
\end{equation}
Corresponding to above computed symmetries, the infinitesimal generator and commutation between them is given as;
\begin{equation}
G_{1}=\frac{\partial }{\partial t}, \quad G_{2}=\frac{\partial }{\partial x},
\end{equation}
\big[ $G_1$,  $G_1$\big]=\big[ $G_2$,  $G_2$\big]=0,  \newline
\big[ $G_1$,  $G_2$\big]=\big[ $G_2$,  $G_1$\big]=0. \newline
\par
Proceeding further, for second case which is considering other terms like P=0, and Q non-zero, symmetries related to this case are
\begin{equation}\begin{aligned}&
\psi=2C_{1}t+C_{o},   \quad \chi=\frac{2C_{1}}{m}=Q(t,x), \\&
\omega=C_{1}x+C_{2}.
\end{aligned}
\end{equation}
other cases could be explored, but considering two cases would be enough because other cases would generate symmetries almost given by the above two cases, we will proceed further for the calculations of commutators, adjoint representations and derivation of optimal system with second case, which one is worth to consider.
\section{lie algebra, commutators and adjoint representations}
Corresponding to above three symmetries of the Chaffee-Infantee equation. We would find lie-algebra and compute the commutators of the given equation, and find the adjoint representations to compute the optimal system of the Chaffee-Infantee equation. Further, in next section, the infinitesimal generators corresponding to the above calculated symmetries for the second case are give as, the three finite dimensional whereas one infinite dimensional group is,
\begin{equation}\begin{aligned}&
G_{1}=\frac{\partial}{\partial x}, \quad \quad G_{2}=\frac{\partial}{\partial t},\\&
G_{3}=2t\frac{\partial}{\partial t}+x\frac{\partial}{\partial x}. \end{aligned}
\end{equation}
The first two symmetries are simply translation in space and time coordinates, while, the infinite dimensional group-generator is,
\begin{equation}
G_{q}=q(t,x)\frac{\partial}{\partial u}.
\end{equation}
The above three infinitesimal generators  form the three-dimensional lie algebra, which is closed, and could be checked by their commutations of the commutators, which is given as; \newline \newline
\big[ $G_1$,  $G_1$\big]=\big[ $G_2$,  $G_2$\big],\big[ $G_2$,  $G_2$\big]=0 \newline
\big[ $G_1$,  $G_2$\big]=\big[ $G_2$,  $G_1$\big]=0, \newline
\big[ $G_1$,  $G_3$\big]=$G_1$, \big[ $G_3$,  $G_1$\big]=$G_1$, \newline
\big[ $G_2$,  $G_3$\big]=$2G_2$, \big[ $G_3$,  $G_2$\big]=$-2G_2$. \newline

\par
Which is closed under Lie-algebra operation. Moving to compute the adjoint representations, which later on would be helpful for computing the optimal system of Chaffee-Infantee equation, the adjoint representations which are computed by the conventional expansion of lie series is,
 \newline \newline
\big[ $G_1$,  $G_1$\big]=$G_1$,\big[ $G_2$,  $G_2$\big]=$G_2$,\big[ $G_2$,  $G_2$\big]=$G_3$\newline
\big[ $G_1$,  $G_2$\big]=$G_2$,\big[ $G_2$,  $G_1$\big]=$G_1$, \newline
\big[ $G_1$,  $G_3$\big]=$G_3-\epsilon G_1$, \big[ $G_3$,  $G_1$\big]=$e^{\epsilon}G_1$, \newline
\big[ $G_2$,  $G_3$\big]=$2G_3-2\epsilon G_2$, \big[ $G_3$,  $G_2$\big]=$e^{2\epsilon }G_2$. \newline

\section{optimal system of Chaffee-Infante equation}
An arbitrary operator for optimal system\cite{PJO} of three dimensional algebra is given as,
\begin{equation}
G=l^1G_{1}+l^2G_{2}+l^3G_{3}.
\end{equation}
Where three constants are choiced, which are arbitrary constants concerning to operators.
\par
The linear transformations corresponding to linear generators, as we will use the form from the generators which after calculations becomes;
\begin{equation}\begin{aligned}&
E_{1}=C^{1}_{13}l^{3}\frac{\partial }{\partial l^1}, \quad E_{2}=C^{2}_{23}l^{3}\frac{\partial }{\partial l^2}, \quad E_{3}=C^{1}_{31}l^{1}\frac{\partial }{\partial l^1}+C^{2}_{32}l^{2}\frac{\partial }{\partial l^2}.
\end{aligned}
\end{equation}
which after calculating coefficients,
\begin{equation*}
C^{1}_{13}=1,C^{2}_{23}=2,C^{1}_{31}=-1,C^{2}_{32}l^{2}=-2,
\end{equation*}
the above equations becomes,
\begin{equation}\begin{aligned}&
E_{1}=l^{3}\frac{\partial }{\partial l^1}, \quad E_{2}=2l^{3}\frac{\partial }{\partial l^2}, \quad E_{3}=(-1)l^{1}\frac{\partial }{\partial l^1}+(-2)l^{2}\frac{\partial }{\partial l^2}.
\end{aligned}
\end{equation}
now we construct the invariant lie equations, and then check the presence of one functinally invariant, after that we will construct optimal system by judicious approach, the lie equations corresponding to three linear transformation equations are,
\begin{equation}
for\quad E_{1}; \quad \hat l^1=l^{3}a_{1}+l^{1}, \hat l^{2}=l^{2},\quad \hat l^{3}=l^{3},
\end{equation}
\begin{equation}
for \quad E_{2}; \quad \hat l^{1}=l^{1}, \hat l^2=2l^{3}a_{2}+l^{2}, \hat l^{3}=l^{3},
\end{equation}
\begin{equation}
for \quad E_{3}; \quad \hat l^{1}=l^{1}a^{-1}_{3}; \hat l^{2}=l^{2}a^{-2}_{3}; \hat l^{3}=l^{3}.
\end{equation}
and the corresponding optimal system by the judicious guess would be
\begin{equation}
G_{1}, \quad G_{2}, \quad kG_{1}+mG_{2}.
\end{equation}
\section{discussion on different types of solutions of Chaffee-Infante equation}
The basic idea of symmetry analysis is to discuss the solutions of the higher non-linear partial differential equations which are otherwise very tough to get, but here after constructing solutions from known solutions, then constructing self-similar or invariant solutions, fundamental solutions or the elementary solutions and infinite ideal solutions corresponding to infinite dimensional algebra is quite important to consider. we transform known solutions to other solutions which are different to the original solutions, being the solutions of the same equation. Similar is the matter with the invariant solutions but same geometry would follow up. here we construct groups and solutions,
As our lie algebra is closed, by the linear combination of given infinitesimal generators\cite{BT} which are
\begin{equation}
C_{1}G_{1}+C_{2}G_{2}+C_{3}G_{3}.
\end{equation}
Choosing different values of constants, dividing it into various cases, we get different one-parameter seven groups which are
\begin{equation}\begin{aligned}&
\Xi_{1}: (x,t,u)\rightarrow (x+\epsilon, t, u), \\&
\Xi_{2}: (x,t,u)\rightarrow (x,t+\epsilon, u),   \\&
\Xi_{3}: (x,t,u)\rightarrow (e^{\epsilon}x, e^{2\epsilon}t ,u), \\&
\Xi_{4}: (x,t,u)\rightarrow (x+\epsilon,t+\epsilon k , u), \\&
\Xi_{5}: (x,t,u)\rightarrow (e^{\epsilon }(1+x)-1, e^{2\epsilon}t, u), \\&
\Xi_{6}: (x,t,u)\rightarrow (e^{\epsilon }(1+x)-1, \frac{1}{2}(e^{2\epsilon}(1+2t)-1), u), \\&
\Xi_{7}: (x,t,u)\rightarrow (e^{\epsilon}x, \frac{1}{2}(e^{2\epsilon}(1+2t)-1), u),  \\&
\Xi_{q}: (x,t,u)\rightarrow (x,t,u+q(t,x)).
\end{aligned}
\end{equation}
Let F(t,x) is the know solution of the given Chaffee-Infantee equation, corresponding to aforementioned one-parameter groups, we can construct solutions for above seven one-parameter and one infinite dimensional group.
\begin{equation}\begin{aligned}&
\Xi_{1}: u=F(x-\epsilon, t, u), \\&
\Xi_{2}: u=F(x, t-\epsilon, u), \\&
\Xi_{3}: u=F(xe^{-\epsilon}, te^{-2\epsilon}, u), \\&
\Xi_{4}: u=F(x-\epsilon, t-\epsilon k, u), \\&
\Xi_{5}: u=F(e^{-\epsilon}(1+x)-1, e^{-2\epsilon}t, u), \\&
\Xi_{6}: u=F(e^{-\epsilon}(1+x)-1, \frac{1}{2}(e^{-2\epsilon}(1+2t)-1, u), \\&
\Xi_{7}: u=F(e^{-\epsilon}x, \frac{1}{2}(e^{-2\epsilon}(1+2t)-1, u), \\&
\Xi_{q}: u=F(t,x)+q(t,x).
\end{aligned}
\end{equation}
Let us construct new solutions from the know solutions of Chaffee-Infante equation\cite{CIE}, the one known soliton  solution of Chaffee-Infante equation is,
\begin{equation}\begin{aligned}&
u(x,t)=-\frac{1}{2}\Big(1+\text{tanh}(\frac{\delta}{2}x+\frac{3\lambda}{4}t)\Big), \\&
 \quad \delta=\sqrt{\frac{\lambda}{2}}. \end{aligned}
\end{equation}
To construct new solutions from known solutions by availing first group, corresponding to first group the transformed solution would be
\begin{equation}
u=-\frac{1}{2}\Big(1+\text{tanh}(\frac{\lambda}{2\sqrt{2}}xe^{-\epsilon}+\frac{3\lambda}{2}te^{-2\epsilon})\Big).
\end{equation}
As epsilon is very small parameter, choosing here epsilon equals to 1, and corresponding to group five, our new transformed solution would be further transformed as,
\begin{equation}
u=-\frac{1}{2}\Big(1+\text{tanh}(\frac{\lambda}{2\sqrt{2}}e^{-1}(e^{-\epsilon }(1+x)-1+\frac{3\lambda}{2}e^{-2}(e^{-2\epsilon}t))\Big ).
\end{equation}
As other new transformed solutions could be generated from aforementioned infinitesimal groups, but the new solutions would be complex in form, though the obtained transformed solutions are the solutions of the original equation but complexity arises, similarly from other groups one can transform further the transformed solutions or one can obtain self-similar or invariant solutions.
\section{reduction of pde's to ode's and power series solutions}
 Employing optimal system for invariant solutions and using infinitesimal groups we would find solutions of partial differential equations, after converting them into ordinary differential equations\cite{GHW}, The exact analytic solutions are found by power series method. \cite{INS}
\subsection{}
Our group invariants determined by characteristic equation of first group, and invariant solution would be
\begin{equation}
t=\eta, \quad u=\mu,
\end{equation}
and
\begin{equation}
 u=f(t), \quad \mu=f(\eta).
\end{equation}
The transformed-equation would be an ODE of first order in time variable 't', Its solution will be found by power series method,
\begin{equation}
\frac{df}{d\eta}+\lambda(f^3-f)=0,
\end{equation}
we let the power series solution of the form,
\begin{equation}
f(\eta)=\sum_{n=0}^{\infty} c_{n}\eta^{n}
\end{equation}
Putting it in above equation gives,
\begin{equation}\begin{aligned}&
c_{1}+\sum_{n=1}^{\infty} c_{n+1}(n+1)\eta^{n}+\lambda \Big(c_{o}^{3}+\sum_{n=1}^{\infty}(\sum_{k=0}^{n}\sum_{j=0}^{k}c_{j}c_{k-j}c_{n-k}\eta^{n})\Big)\\
&-\lambda(c_{o}+\sum_{n=1}^{\infty}c_{n}\eta^{n}).\end{aligned}
\end{equation}
For n=0, the corresponding coefficient would be,
\begin{equation}
c_{1}=\lambda c_{o}-\lambda c_{o}^{3}.
\end{equation}
For $n\geq 1$,
\begin{equation}
c_{n+1}=\frac {-\lambda}{n+1}\Big(c_{n}-(\sum_{k=0}^{n}\sum_{j=0}^{k}c_{j}c_{k-j}c_{n-k})\Big).
\end{equation}
Other coefficients would be determined by above equation, the solution becomes,
\begin{equation}
f(\eta)=c_{o}+c_{1}\eta+\sum_{n=1}^{\infty} c_{n+1}\eta^{n+1}.
\end{equation}
Now replacing into the variable of our target equation,
\begin{equation}
u(x,t)=c_{o}+c_{1}t+\sum_{n=1}^{\infty} c_{n+1}t^{n+1}.
\end{equation}
Putting the values of coefficients, the solution of PDE, is,
\begin{equation}
u(x,t)=c_{o}+(\lambda c_{o}+\lambda c_{o}^{3})t+\sum_{n=1}^{\infty}\frac {-\lambda}{n+1}\Big(c_{n}-(\sum_{k=0}^{n}\sum_{j=0}^{k}c_{j}c_{k-j}c_{n-k})\Big)t^{n+1}.
\end{equation}
Where $c_{o}$ is arbitrary constant.
\subsection{}
Corresponding to second group, we get by solving characteristic equation, the group invariants and invariant solutions
\begin{equation}
x=\eta,\quad u=\mu
\end{equation}
and
\begin{equation}
u=f(x), \quad \mu=f(\eta).
\end{equation}
The given PDE  would be reduced to second order ODE in variable 'x'.
\begin{equation}
\frac{d^2f}{d\eta^2}+\lambda(f^{3}-f)=0,
\end{equation}
Its solution will be found by power series method, we let the power series solution of the form
\begin{equation}
f(\eta)=\sum_{n=0}^{\infty} c_{n}\eta^{n}
\end{equation}
which after substitution in above equation gives,
\begin{equation}\begin{aligned}&
2c_{2}+\sum_{n=1}^{\infty} c_{n+2}(n+1)(n+2)\eta^{n}-\lambda c_{o}^{3}-\lambda \sum_{n=1}^{\infty}(\sum_{k=0}^{n}\sum_{j=0}^{k}c_{j}c_{k-j}c_{n-k}\eta^{n})+\\
&\lambda c_{o}+\lambda \sum_{n=1}^{\infty}c_{n}\eta^{n}.\end{aligned}
\end{equation}
For n=0, the corresponding coefficient would be
\begin{equation}
c_{2}=\frac{\lambda c_{o}^{3}-\lambda c_{o}}{2}.
\end{equation}
For $n\geq 1$,
\begin{equation}
c_{n+2}=\frac {\lambda}{(n+1)(n+2)}(\sum_{k=0}^{n}\sum_{j=0}^{k}c_{j}c_{k-j}c_{n-k}-c_{n}).
\end{equation}
Other coefficients would be determined by above equation, the solution becomes,
\begin{equation}
f(\eta)=c_{o}+c_{1}\eta+c_{2}\eta^{2}+\sum_{n=1}^{\infty} c_{n+2}\eta^{n+2}.
\end{equation}
replacing into the variable of our target equation,
\begin{equation}
u(x,t)=c_{o}+c_{1}x+c_{1}x^{2}+\sum_{n=1}^{\infty} c_{n+2}x^{n+2}
\end{equation}
putting the values of coefficients, the solution of PDE, is,
\begin{equation}\begin{aligned}&
u(x,t)=c_{o}+c_{1}x+\frac {(\lambda c_{o}^{3}-\lambda c_{o})}{2}x+\\
&\sum_{n=1}^{\infty}\frac {\lambda}{(n+1)(n+2)}\Big((\sum_{k=0}^{n}\sum_{j=0}^{k}c_{j}c_{k-j}c_{n-k})-c_{n}\Big)x^{n+2}\end{aligned}
\end{equation}
where $c_{o}$ and $c_{1}$ are arbitrary constants.
\subsection{}
Now, for the third group, by solving characteristic equations, the group invariants are,
\begin{equation}
\frac{t}{x^{2}}=\eta,\quad u=\mu
\end{equation}
and
\begin{equation}
u=f(\frac{t}{x^{2}}), \quad \mu=f(\eta).
\end{equation}
which by substitution reduce to standard second order ODE in $\eta$, The ODE is given as,
\begin{equation}
\frac{1+2x}{x^2}\frac{df}{d\eta}-\frac{4t^2}{x^6}\frac{d^2f}{d\eta^2}+\lambda f^3-\lambda f=0
\end{equation}
Its solution will be found by power series method, we let the power series solution of the form
\begin{equation}
f(\eta)=\sum_{n=0}^{\infty} c_{n}\eta^{n},
\end{equation}
substituting in above equation gives,
\begin{equation}\begin{aligned}&
\frac{1+2x}{x^{2}}c_{1}+\frac{1+2x}{x^{2}}\sum_{n=1}^{\infty} c_{n+1}(n+1)\eta^{n}-\frac{8t^{2}}{x^{6}}c_{2}\\
&-\frac{4t^{2}}{x^{6}}\sum_{n=1}^{\infty} c_{n+2}(n+1)(n+2)\eta^{n}+\lambda c_{o}^{3}+\lambda \sum_{n=1}^{\infty}(\sum_{k=0}^{n}\sum_{j=0}^{k}c_{j}c_{k-j}c_{n-k}\eta^{n})-\\
&\lambda c_{o}-\lambda \sum_{n=1}^{\infty}c_{n}\eta^{n}.\end{aligned}.
\end{equation}
For n=0, the corresponding coefficient would be,
\begin{equation}
c_{2}=\frac{(1+2x)}{8t^{2}}x^{4}c_{1}+x^{6}\frac{(\lambda c_{o}^{3}-\lambda c_{o})}{8t^2}.
\end{equation}
For $n\geq 1$,
\begin{equation}\begin{aligned}&
c_{n+2}=-\frac {x^{6}}{4t^{2}(n+1)(n+2)}(\frac{1+2x}{x^{2}})(n+1)c_{n+1}+\\
&\lambda\Big(\sum_{k=0}^{n}\sum_{j=0}^{k}c_{j}c_{k-j}c_{n-k}-c_{n}\Big).\end{aligned}
\end{equation}
Other coefficients would be determined by above equation, the solution becomes
\begin{equation}
f(\eta)=c_{o}+c_{1}\eta+c_{2}\eta^{2}+\sum_{n=1}^{\infty} c_{n+2}\eta^{n+2}
\end{equation}
replacing into the variable of our target equation,
\begin{equation}
u(x,t)=c_{o}+c_{1}(\frac{t}{x^{2}})+c_{2}(\frac{t}{x^{2}})^{2}+\sum_{n=1}^{\infty} c_{n+2}(\frac{t}{x^{2}})^{n+2}
\end{equation}
putting  the coefficients in above equation, the solution of PDE, is,
\begin{equation}\begin{aligned}&
u(x,t)=c_{o}+c_{1}(\frac{t}{x^{2}})+\frac{(1+2x)}{8t^{2}}x^{4}c_{1}+x^{6}\frac{(\lambda c_{o}^{3}-\lambda c_{o})}{8t^2}(\frac{t}{x^{2}})^{2}\\
&+-\frac {x^{6}}{4t^{2}(n+1)(n+2)}(\frac{1+2x}{x^{2}})(n+1)c_{n+1}+\\
&\lambda\Big(\sum_{k=0}^{n}\sum_{j=0}^{k}c_{j}c_{k-j}c_{n-k}-c_{n}\Big)(\frac{t}{x^{2}})^{n+2}\end{aligned}
\end{equation}
where $c_{o}$ and $c_{1}$ are arbitrary constants.\par
To solve and get solutions of other PDE's, their reduced forms of ODE's are given, which could be solved by similar power series method discussed above.
\subsection{}
To attain solutions from fourth infinitesimal group, the group invariants and invariant solution after solving suitable characteristic equations is given as,
\begin{equation}
\eta=t-kx, \quad u=\mu,
\end{equation}
and
\begin{equation}
 \mu=f(t-kx).
\end{equation}
The corresponding reduced ODE, which could be solved by power series method is given as,
\begin{equation}
\frac{df}{d\eta}-k\frac{d^2f}{d\eta^2}+\lambda f^{3}-\lambda f=0
\end{equation}
 \subsection{}
For the fifth infinitesimal group, solving characteristic equation gives,
\begin{equation}
\eta=\frac{1+x}{\sqrt{t}}, \quad u=\mu.
\end{equation}
And the corresponding invariant solution are given as
\begin{equation}
\mu=f(\eta), \quad u=\mu=f(\frac{1+x}{\sqrt{t}}).
\end{equation}
The reduced PDE, in form of ODE, which could be solved by similar above power series method is,
\begin{equation}
-\frac{1+x}{2\sqrt{t^3}}\frac{df}{d\eta}-\frac{1}{t}\frac{d^2 f}{d \eta^2}+\lambda(f^3-f)=0
\end{equation}
\subsection{}
Now, In case of sixth group the invariants, after solving suitable characteristic equations would be
\begin{equation}
\eta=\frac{1+x}{\sqrt{1+2t}}, \quad u=\mu.
\end{equation}
And the corresponding invariant solution are given as
\begin{equation}
\mu=f(\eta), \quad u=\mu=f(\frac{1+x}{\sqrt{1+2t}}).
\end{equation}
The reduced ODE, which could be solved by above power series method is,
\begin{equation}
-\frac{1+x}{\sqrt{(1+2t)^3}}\frac{df}{d\eta}-\frac{1}{1+2t}\frac{d^2 f}{d\eta^2}+\lambda(f^3-f)=0
\end{equation}
\subsection{}
In case of seventh group the invariants, after solving suitable characteristic equations would be,
\begin{equation}
\eta=\frac{x}{\sqrt{1+2t}}, \quad u=\mu.
\end{equation}
And the corresponding invariant solution are given as
\begin{equation}
\mu=f(\eta), \quad u=\mu=f(\frac{x}{\sqrt{1+2t}}).
\end{equation}
The reduced ODE, which could be solve by power series method to get exact analytic solutions is,
\begin{equation}
-\frac{x}{\sqrt{(1+2t)^3}}\frac{df}{d\eta}-\frac{1}{1+2t}\frac{d^2 f}{d\eta^2}+\lambda(f^3-f)=0
\end{equation}
whereas, the last group, which is infinite dimensional subalgebra, basically generator of an infinite group. This infinite group in most cases could not generate invariant solutions. However, its conveniency lies in getting new solutions form known one.
\section{Conservation Laws}
To construct conservation laws for an infinitesimal symmetry(2), of non-linear Chaffee-Infante equation we employ ibragimov's method of conservation laws\cite{ACTTT}. The pith of it, is stated as, for the adjoint equation which preserves all symmetries of Chaffee-Infante equation, which is,
\begin{equation}
H^{*}=(x, u, v, u_{1}, v_{1}, u_{2}, v_{2})=0,
\end{equation}
Whereas,
\begin{equation}
H^{*}=(x, u, v, u_{1}, v_{1}, u_{2}, v_{2})=\frac{\partial L}{\partial u},
\end{equation}
with L as formal lagrangian. and the variational derivative,
\begin{equation}
\frac{\delta }{\delta u}=\frac{\partial }{\partial u}-D_{i}(\frac{\partial L}{\partial u_{i}})+D_{i}D_{j}(\frac{\partial L}{\partial u_{ij}}),
\end{equation}
Then any symmetry would provide conservation law\cite{ACTTTT} with conserved vectors,
\begin{equation}
T^{i}=\xi L+w\big[\frac{\partial L}{\partial u_{i}}-D{j}(\frac{\partial L}{\partial u_{ij}})\big]+D_{j}(w)\frac{\partial L}{\partial u_{ij}}.
\end{equation}
Generally Construction of conservation laws is not limited only for optimal system. there would be conservation law corresponding to every symmetry of given equation or system of equations. In case of our target equation, The construction of conservation vectors is done for an optimal system. To check adjointness\cite{ACTT}, suppose
\begin{equation}
H=u_t-u_{xx}+\lambda(u^3-u),
\end{equation}
The formal lagrangian,
\begin{equation}
L=v(u_t-u_{xx}+\lambda(u^3-u)),
\end{equation}
Employing self-adjointness condition
\begin{equation}
 \frac{\partial L}{\partial u}=0,
\end{equation}
\begin{equation}
\implies \lambda(3u^{2}v-v)-v_{t}-v{xx}=0.
\end{equation}
is an adjoint equation. which is not strictly self adjoint, As v=u, is not an original equation.
Now, put v=h(u) in adjoint equation(64). It gives,
\begin{equation}
3u^2h(u)\lambda-h(u)\lambda-u_{t}h^\prime(u)-h^{\prime\prime}u_{x}^{2}-h^\prime(u)u_{xx}.
\end{equation}
By using self-adjoint definition;
\begin{equation}\begin{aligned}&
-H^{*}+\Lambda H=\lambda h(u)-\lambda 3u^{2}h(u)+u_{t}h^\prime(u)+h^{\prime\prime}(u)u^2_x+h^\prime(u)u_{xx}+\\
&\Lambda(u_t-u_{xx}+\lambda u^3-\lambda u), \end{aligned}
\end{equation}
\begin{equation}
=u_t(\Lambda+h^\prime(u))+u_{xx}(-\Lambda+h^\prime(u))-\Lambda u+\Lambda u^{3}+h^{\prime\prime}(u)u^2_{x}+\lambda h(u)-3u^{2}h(u)\lambda,
\end{equation}
Coefficients of
\begin{equation}
 u_{t}, u_{xx} \implies  \Lambda+h^\prime(u)=0, \quad -\Lambda+h^\prime(u)=0.
\end{equation}
\begin{equation}
\implies \Lambda=h^\prime(u),
\end{equation}
Shows, target equation is not quasi-self adjoint.
Now, put v=h(t,x,u) in adjoint equation(64). It gives, after above similar procedure,
\begin{equation}
\Lambda=h_{u}(t,x,u)=0,
\end{equation}
\begin{equation}
-h(t,x,u)+3u^{2}h(t,x,u)-h_{t}(t,x,u)-h_{xx}(t,x,u)=0.
\end{equation}
which must be satisfy for two independent and one dependent variable. which is not the case in above equation, so, our given equation is not non-linearly self-adjoint equation. Now for conserved vector of conservation equation we employ formula(10), And focus on symmetries containing in optimal system, the corresponding conserved vectors are given as;
\begin{equation}\begin{aligned}&
\text{For} \quad G_{1}=\frac{\partial }{\partial t}, \quad w=-u{t}, \\&
T^{t}=-vu_{xx}+v\Lambda\lambda u^{3}-v\Lambda\lambda u ,\quad \text{and} \\&
T^{x}=-u_{t}v_{x}+u_{tx}v.\end{aligned}
\end{equation}
\begin{equation}\begin{aligned}&
\text{For} \quad G_{2}=\frac{\partial }{\partial x}, \quad w=-u_{x}, \\
&  T^{t}=-u_{x}v, \text{and} \\
&T^{x}=vu_{t}-vu_{xx}+\Lambda \lambda vu^{3}-\Lambda \lambda vu-u_{x}v_{x}+u_{xx}v.\end{aligned}
\end{equation}
\begin{equation}\begin{aligned}&
\text{For}\quad G_{3}=x\frac{\partial }{\partial x}+ 2t\frac{\partial }{\partial t}, \quad w=-xu_{x}-2tu_{t}, \\
&T^{t}=v(-u_{xx}+\Lambda\lambda u^{3}-\Lambda \lambda u-xu_{x}), \quad \text{and}, \\
&T^{x}=xvu_{t}-xvu_{xx}+x\Lambda\lambda vu^{3}-x\Lambda \lambda vu-xu_{x}-2tu_{t}\\&+
v_{xx}+vu_{x}+vxu_{xx}.
\end{aligned}
\end{equation}
\begin{equation}\begin{aligned}&
\text{For} \quad aG_{1}+G_{2}=a\frac{\partial }{\partial x}+\frac{\partial }{\partial t}, \text\quad{and}\quad  w=-au_{x}-u_{t}.\\
&T^{t}=v(-u_{xx}+\Lambda \lambda u^{3}-\Lambda \lambda u-au_{x}), \text{and} \\
&T^{x}=avu_{t}+\Lambda \lambda avu^{3}-a\Lambda \lambda vu-av+{x}u_{x}-u_{t}v_{x}.
\end{aligned}
\end{equation}
\section{Conclusion}
We performed symmetry analysis of non-linear chaffee infante equation of two dimensions, the results are quite useful especially being usefulness of invariant and transformed solutions, which are new but solutions of the same equation, we further enhanced our analysis by constructing optimal system and at the end we discussed some exact solutions corresponding to infinitesimal generators, further calculations are performed for conservation laws\cite{NH}.Which are quite useful from the physical perspectives of partial differential equations.

\end{document}